\documentclass[12pt]{article}
\usepackage{amssymb}

\def\emline#1#2#3#4#5#6{%
       \put(#1,#2){\special{em:moveto}}%
       \put(#4,#5){\special{em:lineto}}}
\def\newpic#1{}

\author{I. G. Korepanov\\
\normalsize Southern Ural State University\\[-0.5ex]
\normalsize 76 Lenin avenue\\[-0.5ex]
\normalsize 454080 Chelyabinsk, Russia\\[-0.5ex]
\normalsize E-mail: kig@susu.ac.ru}
\date{}
\title{$SL(2)$-solution of the pentagon equation and invariants of three-dimensional
manifolds}

\def\be{\begin{equation}}
\def\ee{\end{equation}}

\def\pa#1#2{{\partial#1\over\partial#2}}

\def\card{\mathop{\rm card}\nolimits}
\def\minor{\mathop{\rm minor}\nolimits}

\begin{document}
\maketitle

\begin{abstract}
Building on a classical solution to the pentagon equation, constructed
earlier by the author and E.~V.~Martyushev and related to the flat
geometry invariant under the group $SL(2)$, we construct an algebraic
complex corresponding to a triangulation of a three-manifold.
In case if this complex is acyclic (which is confirmed by examples), 
we use it for constructing a manifold invariant .

\end{abstract}

\section{Introduction}
\label{sec 1}

Recently, acyclic complexes of a new kind were invented which
can be put in correspondence to triangulations of three- and
four-dimensional manifolds \cite{K4I,K4II,K4III}.
Euclidean metric values (such as edge lengths, dihedral angles and Euclidean 
coordinates of vertices) were ascribed to the elements of these
triangulations. The linear spaces entering in the complex were constructed
out of infinitesimal variations of such values. The result was the construction
of invariants of piecewise-linear manifolds which were expressed in terms
of the torsion of the complex together with volumes of simplices of different
dimensions entering in the triangulation.

The basis for constructing these new acyclic complexes consisted in algebraic
relations corresponding in a natural sense to elementary rebuildings of
a manifold triangulation --- Pachner moves, namely the move
$2\to 3$ (two tetrahedra having a common face are replaced by three tetrahedra
having a common edge; we call such relation {\em pentagon equation\/}) in
the three-dimensional case, and the move $3\to 3$ in the four-dimensional
case. Then, the complex was built in such way as to ensure the construction
of the invariant of {\em all\/} Pachner moves.

In the three-dimensional case, the invariant passes successfully a standard
test of distinguishing the lens spaces~\cite{KM2}. The more complicated
four-dimensional case requires further investigations.

In the present paper, we are dealing again with three-manifolds but, instead
of using the solution to pentagon equation related to Euclidean values,
we are using the {\em $SL(2)$-solution} found in paper~\cite{KM1}. We continue
to use this name for it, although a bigger group arises, too, in our constructions
--- the group of area-preserving affine motions of a plane.

Two representations of the manifold~$M$ fundamental group $\pi_1 (M)$ enter
in a natural way in the construction of acyclic complex. To explain this,
consider once again the Euclidean case. In paper~\cite{K1}, we put the
vertices of a triangulation of manifold~$M$ in a three-dimensional Euclidean
space, while in paper~\cite{KM2} we put there the vertices of the triangulation
of {\em universal cover\/} of manifold~$M$. One can say that every vertex
of $M$'s triangulation was multiplied to $\card \bigl( \pi_1 (M) \bigr)$ its
copies, and the transition from one copy to another was determined by the
image of an element of $\pi_1 (M)$ in the group~$E_3$ of Euclidean motions
of the three-dimensional space with respect to some representation $f\colon
\; \pi_1(M) \to E_3$. In addition to representation~$f$, one can consider
other representations of group $\pi_1(M)$, namely in the linear spaces
of {\em differentials} of which our complex (corresponding to the universal
cover) is built. Here the work with nontrivial representations was initiated
in short note~\cite{M}.

Two similar kinds of representations can be considered in the $SL(2)$-case
as well. The experience of studying the Euclidean case shows that it is
exactly the use of nontrivial representations that leads to the most interesting
manifold invariants. Still, we will confine ourselves in this paper to
the simplest case where both representations are trivial. We hope to study the 
invariants corresponding to nontrivial representations in subsequent papers.

The contents of the remaining sections of this paper is as follows: in
section~\ref{sec 2} we write out the solution to pentagon equation from
\cite{KM1} together with some new ideas. In section~\ref{sec 3} we construct
an algebraic complex on this basis. In section~\ref{sec 4} we study the
behavior of the torsion of the complex (assuming its acyclicity) under
the Pachner moves $2\to 3$ and $1\to 4$, and propose a formula for the
manifold invariant based on this study. In section~\ref{sec 5} the invariant
for sphere~$S^3$ and projective space~$\mathbb RP^3$ is computed (the complexes
turning out indeed acyclic). In the final section~\ref{sec 6} we discuss
the results and plans for future research.

\section{$SL(2)$-solution of pentagon equation}
\label{sec 2}

The Pachner move $2\to 3$ is pictured in Figure~\ref{ris1}:
\begin{figure}
\begin{center}
\unitlength=1.00mm
\special{em:linewidth 0.5pt}
\linethickness{0.5pt}
\begin{picture}(121.00,47.00)
\emline{5.00}{30.00}{1}{20.00}{5.00}{2}
\emline{20.00}{5.00}{3}{35.00}{25.00}{4}
\emline{20.00}{5.00}{5}{50.00}{30.00}{6}
\emline{50.00}{30.00}{7}{35.00}{25.00}{8}
\emline{35.00}{25.00}{9}{5.00}{30.00}{10}
\put(4.00,30.00){\makebox(0,0)[rc]{$A$}}
\special{em:linewidth 0.2pt}
\linethickness{0.2pt}
\emline{5.00}{30.00}{11}{7.00}{30.00}{12}
\emline{9.00}{30.00}{13}{11.00}{30.00}{14}
\emline{13.00}{30.00}{15}{15.00}{30.00}{16}
\emline{17.00}{30.00}{17}{19.00}{30.00}{18}
\emline{21.00}{30.00}{19}{23.00}{30.00}{20}
\emline{25.00}{30.00}{21}{27.00}{30.00}{22}
\emline{29.00}{30.00}{23}{31.00}{30.00}{24}
\emline{33.00}{30.00}{25}{35.00}{30.00}{26}
\emline{37.00}{30.00}{27}{39.00}{30.00}{28}
\emline{41.00}{30.00}{29}{43.00}{30.00}{30}
\emline{45.00}{30.00}{31}{47.00}{30.00}{32}
\emline{49.00}{30.00}{97}{50.00}{30.00}{98}
\put(35.00,24.50){\makebox(0,0)[lt]{$B$}}
\put(51.00,30.00){\makebox(0,0)[lc]{$C$}}
\put(20.00,48.00){\makebox(0,0)[cb]{$D$}}
\put(20.00,4.00){\makebox(0,0)[ct]{$E$}}
\special{em:linewidth 0.5pt}
\linethickness{0.5pt}
\put(59.00,30.00){\vector(1,0){7.00}}
\emline{75.00}{30.00}{35}{90.00}{5.00}{36}
\emline{90.00}{5.00}{37}{105.00}{25.00}{38}
\emline{90.00}{5.00}{39}{120.00}{30.00}{40}
\emline{120.00}{30.00}{41}{105.00}{25.00}{42}
\emline{105.00}{25.00}{43}{75.00}{30.00}{44}
\put(74.00,30.00){\makebox(0,0)[rc]{$A$}}
\special{em:linewidth 0.2pt}
\linethickness{0.2pt}
\emline{75.00}{30.00}{45}{77.00}{30.00}{46}
\emline{79.00}{30.00}{47}{81.00}{30.00}{48}
\emline{83.00}{30.00}{49}{85.00}{30.00}{50}
\emline{87.00}{30.00}{51}{89.00}{30.00}{52}
\emline{91.00}{30.00}{53}{93.00}{30.00}{54}
\emline{95.00}{30.00}{55}{97.00}{30.00}{56}
\emline{99.00}{30.00}{57}{101.00}{30.00}{58}
\emline{103.00}{30.00}{59}{105.00}{30.00}{60}
\emline{107.00}{30.00}{61}{109.00}{30.00}{62}
\emline{111.00}{30.00}{63}{113.00}{30.00}{64}
\emline{115.00}{30.00}{65}{117.00}{30.00}{66}
\emline{119.00}{30.00}{67}{120.00}{30.00}{68}
\put(105.00,24.50){\makebox(0,0)[lt]{$B$}}
\put(121.00,30.00){\makebox(0,0)[lc]{$C$}}
\put(90.00,48.00){\makebox(0,0)[cb]{$D$}}
\put(90.00,4.00){\makebox(0,0)[ct]{$E$}}
\emline{90.00}{5.00}{69}{90.00}{7.00}{70}
\emline{90.00}{9.00}{71}{90.00}{11.00}{72}
\emline{90.00}{13.00}{73}{90.00}{15.00}{74}
\emline{90.00}{17.00}{75}{90.00}{19.00}{76}
\emline{90.00}{21.00}{77}{90.00}{23.00}{78}
\emline{90.00}{25.00}{79}{90.00}{27.00}{80}
\emline{90.00}{29.00}{81}{90.00}{31.00}{82}
\emline{90.00}{33.00}{83}{90.00}{35.00}{84}
\emline{90.00}{37.00}{85}{90.00}{39.00}{86}
\emline{90.00}{41.00}{87}{90.00}{43.00}{88}
\emline{90.00}{45.00}{89}{90.00}{47.00}{90}
\special{em:linewidth 0.5pt}
\linethickness{0.5pt}
\emline{75.00}{30.00}{93}{90.00}{47.00}{94}
\emline{90.00}{47.00}{95}{120.00}{30.00}{96}
\emline{105.00}{25.00}{97}{90.00}{47.00}{98}
\emline{20.00}{47.00}{99}{50.00}{30.00}{100}
\emline{35.00}{25.00}{101}{20.00}{47.00}{102}
\emline{20.00}{47.00}{103}{5.00}{30.00}{104}
\end{picture}
\end{center}
\caption{Move $2\to 3$}
\label{ris1}
\end{figure}
adjacent tetrahedra $EABC$ and $ABCD$ which belong to a triangulation of
a three-dimensional oriented manifold are replaced with three tetrahedra $ABED$, 
$BCED$ and $CAED$. We put in correspondence to every oriented edge a real
number, for instance, number $\lambda_{AB}$ to edge $AB$, and assume that
for all edges
\be
\lambda_{BA}=-\lambda_{AB}.
\label{0a}
\ee

If necessary, we can extend the field to which numbers $\lambda$ and related
values belong to the field~$\mathbb C$ of complex numbers. On the other
hand, we are not considering at this moment an interesting question of
what we will get if $\lambda$'s belong to a field of a {\em finite characteristic}.

For every oriented two-dimensional face, we construct a value $S$ which
is the circulation of values $\lambda$; e.g., for face $ABC$, by definition,
\be
S_{ABC}=\lambda_{AB} + \lambda_{BC} + \lambda_{CA}.
\label{0b}
\ee

Now we ascribe numerical values to dihedral angles at the edges of a tetrahedron. 
Consider an oriented tetrahedron $BCED$ and its oriented edge~$ED$. 
We say that the ``value of the dihedral angle'' at $ED$ is, by definition,
\be
\alpha = \frac{1}{2} \frac{S_{BDC} + S_{BEC}}{S_{BDE} S_{CDE}}.
\label{0c} 
\ee
Thus, $\alpha$ changes its sign both when we change the orientation of
the tetrahedron (without changing the orientation of edge $ED$) and when
we change the orientation of edge $ED$ (without changing the orientation of the
tetrahedron).

Introduce similar values for tetrahedra $CAED$ and $ABED$:
\be
\beta = \frac{1}{2} \frac{S_{CDA} + S_{CEA}}{S_{CDE} S_{ADE}},
\label{0d}
\ee
\be
\gamma = \frac{1}{2} \frac{S_{ADB} + S_{AEB}}{S_{ADE} S_{BDE}}.
\label{0e}
\ee

We assume from now on that (if the contrary is not stated explicitly) we
always choose the {\em positive\/} orientation for all tetrahedra belonging
to the manifold triangulation, i.e., the orientation determined by the
order of tetrahedron vertices coincides with the fixed orientation of the
whole manifold. We set by definition
\be
\omega_{ED}=-\omega_{DE} = \alpha + \beta + \gamma.
\label{0f}
\ee

Under the agreement of the previous paragraph, this value is uniquely
determined by the oriented edge~$ED$. Moreover, its definition generalizes
in the obvious way for the case where edge~$ED$ is common for more than
three tetrahedra. Due to the reasons which will soon be clear, we will
call $\omega_{ED}$ the {\em curvature around edge~$ED$}. 

Return now to Figure~\ref{ris1}. The following formula takes place:
\be
S_{ABC} = S_{ADE}S_{BDE}S_{CDE} \pa{\omega_{ED}}{\lambda_{ED}},
\label{k1}
\ee
where all values are calculated at such $\lambda_{ED}$ for which
$\omega_{ED}=0$. The relation~(\ref{k1}) is exactly the solution of pentagon
equation that lies in the foundation of our subsequent constructions.
In its left-hand side, there is the value $S_{ABC}$ corresponding exactly
to the face which is present in the left diagram of Figure~\ref{ris1} but
absent from the right diagram, and a similar statement is true for the
right-hand side of~(\ref{k1}) as well.

The validity of relation~(\ref{k1}) can be verified directly. 
The ideas that can lead to it can be found in paper~\cite{KM1} (cf.~formula (15) 
of that paper).

The condition $\omega_{ED}=0$ is equivalent to the bilinear relation
\be
S_{ADB} S_{CDE} + S_{BDC} S_{ADE} + S_{CDA} S_{BDE} = 0,
\label{k2}
\ee
and~(\ref{k2}) is transformed into itself under any transposition of letters 
$A,\ldots,E$. On the other hand, the relation~(\ref{k2}) holds if we assume
that $A,\ldots,E$ are points in the usual {\em plane}~$\mathbb R^2$, while
$S_{\ldots}$ are the oriented areas of triangles. Besides the bilinear relations,
the areas obey, of course, linear relations of type $S_{ABC} + S_{ACD} = S_{ABD} +
S_{BCD}$, which is in accordance with such formulas as~(\ref{0b}).

This motivates the following constructions. Let $A$, $B$, $D$ and $E$ be
points in the plane $\mathbb R^2$. We remark at once that we don't need
to calculate, e.g., distances between points or angles (in the Euclidean
sense) within this plane, but only {\em areas\/} of figures. Let all areas
of triangles with vertices in these points be given. Then one can find
that
\be
\overrightarrow{EB} = \frac{S_{EBA} \overrightarrow{ED} + S_{EDB} 
\overrightarrow{EA}}{S_{EDA}},
\label{k3}
\ee
where $S_{EBA}$ is the oriented area of triangle $EBA$ and so~on.

Let us add one more point, $C$, to our four points. Replacing in~(\ref{k3})
the pair $A,B$ first with $B,C$ and then with $C,A$, we get two more relations:
\be
\overrightarrow{EC} = \frac{S_{ECB} \overrightarrow{ED} + S_{EDC} 
\overrightarrow{EB}}{S_{EDB}},
\label{k4}
\ee
\be
\overrightarrow{EA} = \frac{S_{EAC} \overrightarrow{ED} + S_{EDA} 
\overrightarrow{EC}}{S_{EDC}}.
\label{k5}
\ee

Now we look at formulas (\ref{k3}), (\ref{k4}) and (\ref{k5}) from a
different point of view. We forget for a while that $A,\ldots,E$ are
points from $\mathbb R^2$ and regard these letters simply as subscripts
numbering some values $\lambda$ on which only the antisymmetry
condition~(\ref{0a}) is imposed and from which values $S{\ldots}$ are
constructed according to formulas of type~(\ref{0b}).
So, we consider formulas (\ref{k3}),
(\ref{k4}) and (\ref{k5}) simply as linear relations imposed on some nonzero
vectors $\overrightarrow{EA}$, $\overrightarrow{EB}$, $\overrightarrow{EC}$ and
$\overrightarrow{ED}$ lying in~$\mathbb R^2$.

If~(\ref{k2}) holds, then these relation agree with each other. We will
not, however, require the validity of (\ref{k2}), but instead of this we
change $\overrightarrow{EA}$ to $\overrightarrow{EA}_{\hbox{\scriptsize new}}$
in the l.h.s.\ of~(\ref{k5}), and express 
$\overrightarrow{EA}_{\hbox{\scriptsize new}}$ in terms of
$\overrightarrow{ED}$ and $\overrightarrow{EA}$, substituting expressions
(\ref{k4}) and (\ref{k3}) for $\overrightarrow{EC}$ and $\overrightarrow{EB}$.  
We get, after some calculation:
$$
\overrightarrow{EA}_{\hbox{\scriptsize new}} = \overrightarrow{EA} +
\omega_{ED} S_{EDA} \overrightarrow{ED}.
$$

 From the viewpoint of Figure~\ref{ris1}, the relation (\ref{k3}) deals
with tetrahedron $ABED$, (\ref{k4}) --- with tetrahedron $BCED$ and (\ref{k5}) ---
with~$CAED$. Their successive use corresponds, one can say, to going around the
edge~$ED$. Namely, we see that the element of $SL(2)$ which corresponds to
this going around is determined by the following transformation of
(unnormed) bases in the plane~$\mathbb R^2$:
\be
(\overrightarrow{ED},
\overrightarrow{EA}) \to (\overrightarrow{ED}, \overrightarrow{EA}+\omega_{ED}
S_{EDA}\overrightarrow{ED} ).
\label{kak-to formuloj}
\ee

If we return to considering $E$, $D$ and~$A$ in formula~(\ref{kak-to formuloj}) as
points in the plane~$\mathbb R^2$, $\overrightarrow{ED}$ and~$\overrightarrow{EA}$
--- as vectors with corresponding origins and ends, and $S_{EDA}$ --- as
the area of triangle $EDA$, then the transformation~(\ref{kak-to formuloj})
depends {\em only on the vector $\overrightarrow{ED}$ and number $\omega_{ED}$},
and not on the vector~$\overrightarrow{EA}$.

\section{Construction of the algebraic complex}
\label{sec 3}

The algebraic complex that we are going to build in this Section has the
following form:
\be
0\to \mathbb R^6 \stackrel{f_{1}}{\to} \mathbb R^{\hbox{\scriptsize $3\#$vertices}} 
\stackrel{f_{2}}{\to} \mathbb R^{\hbox{\scriptsize $\#$edges}} 
\stackrel{f_{3}}{\to} \mathbb R^{\hbox{\scriptsize $\#$edges}} 
\stackrel{f_{4}}{\to} \mathbb R^{\hbox{\scriptsize $3\#$vertices}} 
\stackrel{f_{5}}{\to} \mathbb R^6 \to 0.
\label{l1}
\ee

Here $\#$vertices is, of course, the number of vertices in the simplicial
complex, while $\#$edges --- the number of edges. The bases are fixed in
all spaces, so the mappings are identified with matrices.
We now pass on to defining the mappings $f_i$ in the sequence~(\ref{l1}).

We put all the vertices of the complex in the two-dimensional coordinate
space~$\mathbb R^2$, i.e.\ we ascribe to every vertex~$A$ real numbers
$x_A$ and~$y_A$. We assume that values $x_A$ and~$y_A$ can take
infinitesimal variations $dx_A$ and~$dy_A$. Besides, we put in correspondence
to every vertex~$A$ one more real number~$\varkappa_A$ which, too, can
take a variation~$d\varkappa_A$. We understand the {\em second left\/}
nonzero space $\mathbb R^{\hbox{\scriptsize $3\#$vertices}}$ in the sequence~(\ref{l1})
as consisting of columns of differentials $(dx_A,\allowbreak dy_A,\allowbreak 
d\varkappa_A,\allowbreak dx_B,\allowbreak dy_B,\allowbreak d\varkappa_B,\allowbreak 
\ldots)^{\rm T}$ (the superscript $\rm T$ means matrix transposing). 

The {\em leftmost\/} nonzero space $\mathbb R^6$ is the five-dimensional
Lie algebra of infinitesimal affine area-preserving transformations of the 
plane $\mathbb R^2$ plus (direct sum) the one-dimensional space of 
differentials~$d\varkappa$. A vector in this space is represented by a
column
\be
(dt_1, dt_2, dt_3, dx, dy, d\varkappa)^{\rm T},
\label{l2}
\ee
which corresponds to the element $\pmatrix{dt_3 & dt_1 \cr dt_2 & -dt_3}$
of Lie algebra~$\mathfrak s\mathfrak l(2)$, the translation by vector 
$\pmatrix{dx \cr dy}$ in $\mathbb R^2$ and the translation of parameter
$\varkappa$ by $d\varkappa$. The mapping~$f_1$ in sequence~(\ref{l1}) is
as follows: to a vector~(\ref{l2}) there correspond, in every vertex~$A$,
the vector
\be
\pmatrix{dx_A \cr dy_A} = 
\pmatrix{dt_3 & dt_1 \cr dt_2 & -dt_3} \pmatrix{x_A \cr y_A} +
\pmatrix{dx \cr dy}
\label{m1}
\ee
and
\be
d\varkappa_A = d\varkappa + \frac{x_A\, dy -y_A\, dx}{2}.
\label{m2}
\ee
The meaning of formulas (\ref{m1}) and (\ref{m2}) will become clear when
we start proving that the sequence (\ref{l1}) is a complex.

Now we ascribe to every edge $AB$ the value
\be
\lambda_{AB} = S_{OAB} + \varkappa_B - \varkappa_A,
\label{l3}
\ee
where $S_{OAB}$ is the (oriented) area of triangle $OAB$, $O$~is the origin
of coordinates. The left one of the two spaces $\mathbb R^{\hbox{\scriptsize $\#$edges}}$
in (\ref{l1}) consists of columns of differentials $d\lambda_{AB}$ for
all edges. The mapping~$f_2$, by definition, is obtained by differentiating
the equality~(\ref{l3}), i.e.\ the value $d\lambda_{AB}$ which is obtained
from given $dx_{\ldots},dy_{\ldots}$ and $d\varkappa_{\ldots}$ by means
of $f_2$ is
\be
d\lambda_{AB} = \frac{y_B}{2} dx_A - \frac{x_B}{2} dy_A - \frac{y_A}{2}
dx_B + \frac{x_A}{2} dy_B + d\varkappa_B - d\varkappa_A.
\label{m3}
\ee

As the discussion in Section~\ref{sec 2} shows, the curvatures $\omega_{AB}$
around all edges are zero if all $\lambda_{AB}$ are obtained from vertex
coordinates by formulas~(\ref{l3}) (the adding of values $\varkappa_{\ldots}$ 
obviously does not interfere with the vanishing of the curvatures).
Now we introduce into consideration the differentials of all curvatures.
They will form the {\em right\/} one of the two spaces 
$\mathbb R^{\hbox{\scriptsize $\#$edges}}$ in sequence~(\ref{l1}).
The matrix of mapping $f_3$ is, by definition, $(\partial \omega_a / \partial
\lambda_b)$, where $a$ numbers the edges of the complex and at the same
time the rows of this matrix, while $b$ --- the edges and the columns.
The partial derivatives are calculated on the basis of relation (\ref{0f})
or its generalization to greater number of tetrahedra, the ``angles'' $\alpha$,
$\beta$, $\gamma,\ldots$ being calculated according to formulas of type
(\ref{0c}), (\ref{0d}), (\ref{0e}), where, of course, values $\lambda_{\ldots}$
are substituted according to formulas (\ref{0b}) and~(\ref{0a}).

The right one of the two spaces $\mathbb R^{\hbox{\scriptsize $3\#$vertices}}$ in
sequence~(\ref{l1}) will be the direct sum of Lie algebras $\mathfrak s
\mathfrak l (2)$, one copy of the algebra for every vertex in the complex. 
The mapping $f_4$ is constructed in the following way. Consider, for instance,
a vertex~$E$. The element $d\gamma_E\in \mathfrak s \mathfrak l (2)$ in
which $f_4$ transforms a given configuration of values $d\omega$ on the
edges of complex equals, by definition, the sum over all edges {\em beginning\/} 
in $E$ of the elements of this Lie algebra corresponding to infinitesimal
basis transformations of type~(\ref{kak-to formuloj}).

Assume we are considering an edge $ED$ whose coordinates are
\be
\pmatrix{x_{ED} \cr y_{ED}} = 
\pmatrix{x_D - x_E \cr y_D - y_E}.
\label{l4}
\ee
After some calculation, we find from formula (\ref{kak-to formuloj}), where
we replace $\omega_{ED}$ with $d\omega_{ED}$ and consider, instead of the
transformation from the {\em group}~$SL(2)$, its differential (i.e., subtract
the identity matrix), the following element of algebra
$\mathfrak s \mathfrak l (2)$:
\be
\frac{1}{2} d\omega_{ED}
\pmatrix{-x_{ED}y_{ED} & x_{ED}^2 \cr -y_{ED}^2 & x_{ED}y_{ED}}.
\label{b1}
\ee
To construct the algebraic complex, we need only the structure of vector
space, thus, we will identify the expression~(\ref{b1}) with the column
vector
\be
d\gamma_{ED} = \frac{d\omega_{ED}}{2} \pmatrix{ x_{ED}^2 \cr  x_{ED}y_{ED}
\cr  y_{ED}^2}.
\label{b2}
\ee
Summing over all vertices $D$ joined with $E$ by the edges, we get the
three components of vector
$d\gamma\in \mathbb R^{\hbox{\scriptsize $3\#$vertices}}$ corresponding
to vertex~$E$:
\be
(d\gamma_E)_1 = \sum_D x_{ED}^2 \frac{d\omega_{ED}}{2},
\label{b3}
\ee
\be
(d\gamma_E)_2 = \sum_D x_{ED} y_{ED} \frac{d\omega_{ED}}{2},
\label{b4}
\ee
\be
(d\gamma_E)_3 = \sum_D y_{ED}^2 \frac{d\omega_{ED}}{2}.
\label{b5}
\ee

It remains to construct the right space $\mathbb R^6$ and mapping~$f_5$.
The space will consist of column vectors with components $d\beta_1,\ldots,
\allowbreak d\beta_6$.
By definition, values $d\beta$ obtained by means of $f_5$ from given
$d\gamma$ are (the sums are taken over all vertices $A$ in the complex):
\be
d\beta_1 = \sum_A (d\gamma_A)_1, \quad
d\beta_2 = \sum_A (d\gamma_A)_2, \quad
d\beta_3 = \sum_A (d\gamma_A)_3,
\label{b6}
\ee
\be
d\beta_4 = \sum_A \bigl( y_A (d\gamma_A)_1 - x_A (d\gamma_A)_2 \bigr),
\label{b7}
\ee
\be
d\beta_5 = \sum_A \bigl( y_A (d\gamma_A)_2 - x_A (d\gamma_A)_3 \bigr),
\label{b8}
\ee
\be
d\beta_6 = \sum_A \bigl( y_A^2 (d\gamma_A)_1 - 2x_A y_A (d\gamma_A)_2 + 
x_A^2 (d\gamma_A)_3 \bigr).
\label{b9}
\ee

The sequence (\ref{l1}) is constructed. We will also use for the linear
spaces entering in it somewhat looser but convenient notations in the style
of papers \cite{K4II} and~\cite{K4III}, and write it the following way:
\be
0\to \bigl( \mathfrak s \mathfrak l (2) \hbox{ and translations} \bigr)
\stackrel{f_{1}}{\to} \pmatrix{dx_A\cr dy_A\cr d\varkappa_A}
\stackrel{f_{2}}{\to} (d\lambda_{AB})
\stackrel{f_{3}}{\to} (d\omega_{AB})
\stackrel{f_{4}}{\to} (d\gamma_A)
\stackrel{f_{5}}{\to} (d\beta) \to 0.
\label{b10}
\ee
Here the word ``translations'' in the leftmost nonzero space means
``global translations'' $dx$, $dy$ and $d\varkappa$ acting according to
formulas (\ref{m1}) and~(\ref{m2}); the next space consists of differentials $dx_A$,
$dy_A$, $d\varkappa_A$ for every vertex~$A$ and so on.

Now we prove that the sequence (\ref{b10}) is indeed a complex, that is
the composition of any two successive mappings is zero. We start with the
composition $f_2 \circ f_1$. Substitute in (\ref{m3}) expressions (\ref{m1})
and (\ref{m2}) for $dx_A$, $dy_A$ and~$d\varkappa_A$, as well as similar
expressions for $dx_B$, $dy_B$ and~$d\varkappa_B$. The result of a direct calculation
is then $d\lambda_{AB}=0$. We see that the terms added to $d\varkappa$
in the right-hand side of (\ref{m2}) are chosen so as to compensate the
change of the area of triangle $OAB$ arising when its vertices $A$ and~$B$
are shifted by the vector~$\pmatrix{dx\cr dy}$ while the origin of coordinates $O$
remains in its place.

The fact that $f_3 \circ f_2 = 0$ is evident from geometric considerations:
if the areas of all triangles are determined by their coordinates (in other
words, all the triangles can be placed in $\mathbb R^2$), then all the
curvatures are zero. 

The fact that $f_4 \circ f_3 = 0$ is also evident from geometric considerations.
Imagine a vertex~$E$ and edges going out of it. A transformation of the
basis in the plane~$\mathbb R^2$ can be put in correspondence to any closed path that
does not intersect edges. Indeed, as formulas (\ref{k3}) and~(\ref{k4})
show, the system of coordinates in~$\mathbb R^2$ is uniquely extended from
a given tetrahedron to an adjacent one (having a common face with the first one),
but if the curvatures are nonzero we can get a new system of coordinates
on returning to the initial tetrahedron (which is demonstrated by 
formula~(\ref{kak-to formuloj})). If, however, the path can be contracted
into a point in such way that it does not intersect edges at any moment,
we must, of course, arrive at the system of coordinates identical to the
initial one. For the case of infinitesimal curvatures $d\omega$, this means
that the sum of algebra $\mathfrak s \mathfrak l (2)$ elements of the form
(\ref{b1}) over all vertices~$D$ joined by an edge with~$E$ will be zero
for any~$d\lambda$. Thus, the right-hand sides of expressions (\ref{b3}),
(\ref{b4}) and~(\ref{b5}) will be zero, as desired.

It remains to check that $f_5 \circ f_4 = 0$. We substitute the expressions
(\ref{b3})--(\ref{b5}) in formulas (\ref{b6})--(\ref{b9}), take into account
the relations of type (\ref{l4}) and get zero. This is quite evident for
the expressions (\ref{b6}): the curvature at each edge gives, in view of
$d\omega_{AB}=-d\omega_{BA}$, the mutually opposite contributions at its
two ends; on summing according to (\ref{b6}), one gets zero. The explanation
is a bit harder for (\ref{b7}), (\ref{b8}) and (\ref{b9}); we have to say
simply that these expressions were specially invented in such way as to get zeros.

\section{The behavior of the torsion under moves $2\to 3$ and $1\to 4$
and the formula for the invariant}
\label{sec 4}

There are reasons to believe (see the next section of the present
paper) that at least in many interesting cases the complex (\ref{b10})
is acyclic, i.e.\ the images of mappings coincide exactly with the kernels
of the next mappings. If this is true then one can define the {\em
torsion}~$\tau$ of the complex (\ref{b10}) as the product of some minors
in matrices of mappings $f_1,\ldots,f_5$ taken in the alternating powers
$+1$ and~$-1$. Recall that the bases in all spaces are fixed. Some of the
basis vectors in each space correspond to the rows of the minor belonging
to the left (from this space) mapping, while the rest of them --- to the
columns of the minor belonging to the right mapping.

We choose the exponents $+1$ for the mappings with odd numbers and $-1$
for those with even numbers:
\be
\tau = \frac{\minor f_1 \cdot \minor f_3 \cdot \minor f_5}{\minor f_2 \cdot
\minor f_4}.
\label{t1}
\ee

Now we consider the local rebuilding of type $2\to 3$ of the simplicial complex,
pictured in Figure~\ref{ris1}. The new edge $DE$ adds one new basis vector
in both spaces $(d\lambda_{AB})$ and $(d\omega_{AB})$. We include the corresponding
column and row in $\minor f_3$ (thus, the remaining minors in (\ref{t1})
rest intact). Considerations using the triangular form of the matrices
and similar to those given in section~3 of paper~\cite{K1}, show that this makes
$\minor f_3$ to be multiplied by $\partial \omega_{ED} / \partial \lambda_{ED}$. 
Taking into account formula (\ref{k1}), this shows that the value
\be
\tau \cdot \prod_{\hbox{\scriptsize
\begin{tabular}{c}
over all\\ 
two-dimensional faces
\end{tabular}} } S
\label{t2}
\ee
does not change under the move~$2\to 3$.

Note that we define both the torsion and our future invariant 
{\em to within a sign}. This allows us not to care about the orientation
of the two-dimensional faces in the complex.

Consider now the rebuilding of type $1\to 4$: we add the new vertex~$E$
inside the tetrahedron $ABCD$. In the same way as for the move $2\to 3$,
we add to the minor of mapping~$f_3$ a row and a column corresponding to edge~$ED$.
We then add to the minor of mapping~$f_2$ rows corresponding 
to the three remaining edges, $EA$, $EB$ and~$EC$, and columns corresponding
to $dx_E$, $dy_E$ and~$d\varkappa_E$. Similarly, we add to the minor of
mapping~$f_4$ columns corresponding to $EA$, $EB$ and~$EC$ and three rows
corresponding to the three components of~$d\gamma_E$. Again, considerations
using the triangular form of matrices show that $\minor f_2$ and $\minor f_4$
are simply multiplied by $3\times 3$ minors corresponding to the new rows
and columns. The first of these $3\times 3$ minors is calculated using
formula (\ref{m3}) and turns out to equal $\frac{1}{2} S_{ABC}$, while
the second one is calculated using formulas (\ref{b3}), (\ref{b4}) and~(\ref{b5})
and turns out to equal $S_{EAB} S_{EBC} S_{ECA}$ (both --- to within their
signs), which makes the whole torsion (where we must take into account
also the minor of~$f_3$ which behaves the same way as under a move $2\to 3$)
to be multiplied by
$\displaystyle 2 \prod_{\hbox{\scriptsize
\begin{tabular}{c}
over new\\ 
faces
\end{tabular}} } S^{-1}.$

This all shows that the following value is invariant under all Pachner
moves and can be thus attributed to the manifold~$M$ itself:
\be
I_{SL(2)}(M) = \tau \cdot 
\prod_{\hbox{\scriptsize
\begin{tabular}{c}
over all\\ 
faces
\end{tabular}} } S \cdot 2^{\hbox{\scriptsize $-\#$vertices$-1$}}.
\label{t3}
\ee
Minus one is added to the exponent of number 2 in order that the invariant
be equal to 1 for the sphere~$S^3$, see below Subsection~\ref{subsec 5.1}.

\section{Examples}
\label{sec 5}

\subsection{Sphere $S^3$}
\label{subsec 5.1}

We take the triangulation of the sphere consisting of two tetrahedra, each
of which has vertices $A$, $B$, $C$ and~$D$. We must ascribe to these vertices
some coordinates in the plane~$\mathbb R^2$ (on which, of course, our determinant
will not depend); for instance, we can take: $A(0,0)$, $B(x_B,0)$, $C(0,y_C)$
and $D(x_D,y_D)$.

In the same way as in the Euclidean case (see Section~6 of paper~\cite{K1}), the mapping
$f_3\colon \; (d\lambda) \to  (d\omega)$ is the identical zero, because
the curvature~$\omega$ around any of the six edges is made up of two mutually
opposite summands. Thus, $\minor f_3=1$ (as the determinant of a $0\times 0$
matrix). To calculate the minors of $f_1$ and $f_2$, we must choose six
basis differentials in the space of values
$\pmatrix{dx_A\cr dy_A\cr d\varkappa_A}$ for the rows of minor~$f_1$ (the
calculations go very easily if we choose 
$dx_A$, $dy_A$, $d\varkappa_A$, $dx_B$, $dy_B$ and~$dx_C$), and use the
remaining six basis differentials for columns of minor~$f_2$. The result
is:
\be
\frac{\minor f_1}{\minor f_2} = \pm 8.
\label{v1}
\ee

The calculations in the right side of the sequence are similar although
slightly harder. Their result is:
\be
\frac{\minor f_5}{\minor f_4} = \pm \frac{4}{S_{ABC}S_{ABD}S_{ACD}S_{BCD}}.
\label{v2}
\ee
Substituting the right-hand sides of relations (\ref{v1}) and (\ref{v2})
in formula (\ref{t1}) and then in formula (\ref{t3}),
we get the announced result:
$$
I_{SL(2)}(S^3)=1
$$
(here we omit the $\pm$ sign).

\subsection{Projective space $\mathbb RP^3$}
\label{subsec 5.2}

\begin{figure}
\begin{center}
\unitlength=1mm
\special{em:linewidth 0.5pt}
\linethickness{0.5pt}
\begin{picture}(66.00,86.00)
\emline{5.00}{45.00}{1}{25.00}{35.00}{2}
\emline{25.00}{35.00}{3}{65.00}{45.00}{4}
\emline{65.00}{45.00}{5}{35.00}{5.00}{6}
\emline{35.00}{5.00}{7}{25.00}{35.00}{8}
\emline{35.00}{5.00}{9}{5.00}{45.00}{10}
\emline{5.00}{45.00}{11}{35.00}{85.00}{12}
\emline{35.00}{85.00}{13}{25.00}{35.00}{14}
\emline{65.00}{45.00}{15}{35.00}{85.00}{16}
\special{em:linewidth 0.15pt}
\linethickness{0.15pt}
\emline{5.00}{45.00}{17}{45.00}{55.00}{18}
\emline{45.00}{55.00}{19}{65.00}{45.00}{20}
\emline{65.00}{45.00}{21}{5.00}{45.00}{22}
\emline{25.00}{35.00}{23}{45.00}{55.00}{24}
\emline{45.00}{55.00}{25}{35.00}{85.00}{26}
\emline{35.00}{85.00}{27}{35.00}{5.00}{28}
\emline{35.00}{5.00}{29}{45.00}{55.00}{30}
\put(22.00,69.00){\makebox(0,0)[rb]{$g$}}
\put(48.00,69.00){\makebox(0,0)[lb]{$g'$}}
\put(30.50,65.50){\makebox(0,0)[rb]{$f$}}
\put(41.50,66.50){\makebox(0,0)[lb]{$f'$}}
\put(36.00,64.00){\makebox(0,0)[lb]{$d$}}
\put(20.00,23.00){\makebox(0,0)[rt]{$g'$}}
\put(28.50,23.50){\makebox(0,0)[rt]{$f'$}}
\put(34.50,25.00){\makebox(0,0)[rt]{$d'$}}
\put(39.50,25.50){\makebox(0,0)[lt]{$f$}}
\put(49.50,23.00){\makebox(0,0)[lt]{$g$}}
\put(17.00,40.00){\makebox(0,0)[lb]{$h'$}}
\put(49.00,40.00){\makebox(0,0)[lt]{$h$}}
\put(54.50,50.50){\makebox(0,0)[lb]{$h'$}}
\put(21.00,50.00){\makebox(0,0)[rb]{$h$}}
\put(22.50,45.50){\makebox(0,0)[cb]{$b'$}}
\put(50.00,46.00){\makebox(0,0)[cb]{$b$}}
\put(39.00,50.00){\makebox(0,0)[rb]{$c$}}
\put(30.50,40.50){\makebox(0,0)[rb]{$c'$}}
\put(34.00,46.00){\makebox(0,0)[rb]{$A$}}
\put(4.00,45.00){\makebox(0,0)[rc]{$B$}}
\put(66.00,45.00){\makebox(0,0)[lc]{$B$}}
\put(46.00,55.00){\makebox(0,0)[lb]{$C$}}
\put(24.00,34.00){\makebox(0,0)[rt]{$C$}}
\put(35.00,4.00){\makebox(0,0)[ct]{$D$}}
\put(35.00,86.00){\makebox(0,0)[cb]{$D$}}
\end{picture}
\end{center}
\caption{Triangulation of $\mathbb R P^3$}
\label{ris2}
\end{figure}
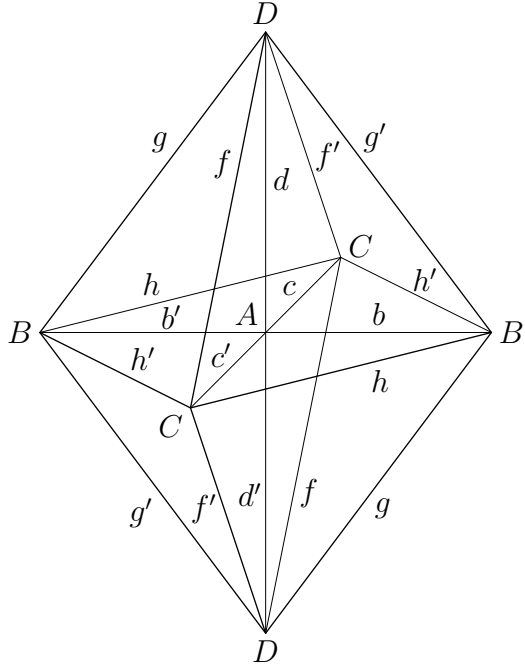

We take the same triangulation as in section 6 of paper~\cite{K1}, where
we considered the Euclidean case.
Namely, there are again 4 vertices $A$, $B$, $C$ and~$D$, but now 12 edges, 16
two-dimensional faces and 8 tetrahedra, see Figure~\ref{ris2}.
Here we have, for example, two edges $b$ and~$b'$ instead of the single edge
$AB$ in Subsection~\ref{subsec 5.1} of the present work, and so on. Note that any two edges
bearing identical notations in Figure~\ref{ris2}, for instance, the two
copies of edge~$f$, are identified, but $f$ and~$f'$ are {\em different\/} edges.

Despite all these differences, formulas (\ref{v1}) and (\ref{v2}) {\em
remain valid\/} for $\mathbb RP^3$ as well. We must only specify that we
employ the six ``primed'' differentials $d\lambda$ for the rows of minor~$f_2$, 
that is, $d\lambda_{b'},\ldots,d\lambda_{f'}$ in place of
$d\lambda_{AB},\ldots,d\lambda_{CD}$ from Subsection~\ref{subsec 5.1}.
Similarly, we use the ``primed'' differentials~$d\omega$ for constructing
the columns of minor~$f_4$.

``Unprimed'' $d\lambda$ and $d\omega$ remain for constructing the minor
of mapping~$f_3$ consisting of the partial derivatives of the six values
$\omega_b$,  $\omega_c$,  $\omega_d$,  $\omega_h$, $\omega_g$ and $\omega_f$
with respect to the variables
$\lambda_b$, $\lambda_c$, $\lambda_d$, $\lambda_h$, $\lambda_g$, $\lambda_f$.
This $6\times 6$ minor is equal to the product of its six elements, namely,
\be
\pa{\omega_b}{\lambda_f},\quad 
\pa{\omega_f}{\lambda_b},\quad 
\pa{\omega_c}{\lambda_g},\quad 
\pa{\omega_g}{\lambda_c},\quad 
\pa{\omega_d}{\lambda_h} \hbox{ \ and \ } 
\pa{\omega_h}{\lambda_d},
\label{v3}
\ee 
while the rest of its elements are zero. The reason for this is that these
remaining elements are made up of pairs of mutually opposite summands corresponding
to tetrahedra having opposite orientations (we are, essentially, repeating
the argumentation from section~6 of paper~\cite{K1}). As for the derivatives
(\ref{v3}), each of them is a sum of two identical terms, corresponding
to two tetrahedra with the same orientation. A derivative of such sort
for one tetrahedron is given in paper~\cite{KM1}, formula~(10). Changing
the notations of that formula to those of the present paper, multiplying
by~2 and ignoring the possible minus sign, we get, for example, the formula 
$$
\pa{\omega_b}{\lambda_f} = \frac{2}{S_{ABC} S_{ABD}}
$$
(it is also not very difficult to deduce this from formulas of type 
(\ref{0c}), (\ref{0d}), (\ref{0e})).

For the whole minor, we get
$$
\minor f_3 = \frac{64}{S_{ABC}^3 S_{ABD}^3 S_{ACD}^3 S_{BCD}^3 }.
$$
Hence,
\be
I_{SL(2)}(\mathbb RP^3) = 64.
\label{v4}
\ee

\section{Discussion}
\label{sec 6}

Our result (\ref{v4}) confirms the hypothesis stated in the end of paper~\cite{KM1}: 
in the case of trivial representations of the fundamental group, the
Euclidean and $SL(2)$ invariants yield the same result (the Euclidean
invariant for~$\mathbb RP^3$, as calculated in~\cite{K1}, was equal to $1/8$,
but this was simply because we defined it there in a slightly different way,
analogous to our present definition in the power~$(-1/2)$).

In the terminology of paper \cite{KM1}, we have globalized the ``local'' pentagon
equation in the present paper, i.e., we have shown how it can be used for
studying not only a local rebuilding $2\to 3$ but the whole manifold. This
globalization turned out to be harder than we expected when writing the
paper~\cite{KM1}. Recall that, in the Euclidean case, the algebraic complex
analogous to~(\ref{b10}) was symmetric with respect to its middle in the
following sense: the matrices of mappings equidistant from the ends of
the complex could be obtained from each other by means of transposing (see the sections
of papers \cite{K4II} and~\cite{K4III} devoted to the three-dimensional
case). For the complex (\ref{b10}), not only this property is no longer
valid, but even its analogue of any kind could not be found as yet. In
general, the feeling is that the geometric and algebraic sense of sequence (\ref{b10})
has not yet been discovered in full.

It seems that the exactness of sequence (\ref{b10}) in some of its terms can be
shown for arbitrary manifolds using ideas parallel to those of section~2 in
paper~\cite{K4II}. There is also another idea: to consider the simplicial
complex corresponding to the {\em universal cover\/} of the manifold, construct
for it the algebraic complex of type (\ref{b10}), and then consider the
subcomplexes of the algebraic complex. They must correspond to different
representations of the fundamental group (which we have discussed in the
Introduction). It looks plausible that this may help both to prove the
acyclicity and construct new manifold invariants. In the Euclidean
case, such activity was initiated in paper~\cite{M}.

Finally, there is a very intriguing question about possible existence of
some quantum relations from which our solution to pentagon equation can
be obtained as a semiclassical limit. This question is motivated by the
fact that, in the Euclidean case, our solution to pentagon equation~\cite{K1}
can be obtained by a limiting procedure from the quantum $6j$-symbols.
Note, however, that even in the Euclidean case the similar question remains
open for manifolds of dimensionality more than three~\cite{K4I,K4II,K4III}.

\medskip

{\bf Acknowledgements. }The work has been performed with the partial financial
support from Russian Foundation for Basic Research, Grant no.~01-01-00059.

\end{document}